\documentclass[12pt]{article}
\usepackage{amsmath,amssymb,amsthm,amscd}
\usepackage[small,nohug,heads=littlevee,noPostScript]{diagrams}
\title{Wreath products,  nilpotent orbits and symplectic deformations}
\author{Baohua Fu}

\newtheorem{Prop}{Proposition}

\newtheorem{Rque}{Remark}

\def\cit{{\mathbb C}}

\def\zit{{\mathbb Z}}
\def\pit{{\mathbb P}}

\def\0{{\mathcal O}}
\def\g{{\mathfrak g}}
\def\h{{\mathfrak h}}
\def\p{{\mathfrak p}}

\def\c{{\mathfrak c}}

\def\Y{{\mathcal Y}}
\def\X{{\mathcal X}}
\def\Z{{\mathcal Z}}

\def\T{{\mathcal T}}

\def\u{{\mathfrak u}}

\def\Hilb{\mathop{\rm Hilb}\nolimits}
\def\Sym{\mathop{\rm Sym}\nolimits}
\def\Sp{\mathop{\rm Sp}\nolimits}
\def\rk{\mathop{\rm rk}\nolimits}
\def\tr{\mathop{\rm tr}\nolimits}
\begin{document}
\maketitle
\begin{abstract}
We recover the wreath product $X:=\Sym^2(\cit^2/{\pm 1})$ as a transversal slice to a nilpotent orbit
in $\mathfrak{sp}_6$. By using deformations of Springer resolutions, we construct a symplectic
deformation of symplectic resolutions of $X$. AMS Classification: 14E15, 14M17
\end{abstract}

{\bf 0. Introduction} \vspace{0.3 cm}
 
Let $H \subset \Sp(2n)$ be a finite sub-group and $X:= \cit^{2n}/H$ the quotient symplectic variety.
Given a projective symplectic resolution 
\begin{equation} \label{reso}
 Z \to X,
\end{equation}
 it was shown in \cite{GK} that
there exists  a symplectic deformation of (\ref{reso}) over $B: = H^2(Z, \cit)$, i.e.  a morphism 
$\Pi: \Z \to \X$ over $B$ such that over the origin $0 \in B$,  $\Pi_0: \Z_0 \to \X_0$ is the resolution (\ref{reso}), and 
over a generic point $b \in B$, $\Z_b, \X_b$ are symplectic smooth varieties
isomorphic under $\Pi_b$, where $\Pi_b$ 
is the restriction of $\Pi$ to the fibers over $b$.
The proof of this theorem is based on the infinitesimal and formal deformations of $\pi$ developed in
\cite{KV} and the globalization is obtained by using the expanding $\cit^*$-action on $X$.
As noted already in \cite{GK}, this deformation is very similar to the deformation of the Springer resolution of nilpotent cones given by 
Grothendieck's simultaneous resolution (\cite{Slo}). However, the construction of symplectic deformations in general is rather implicit.
The purpose of this note is to provide some explicit examples of such deformations.

A class of important examples of symplectic resolutions 
is given by Hilbert-Chow morphisms (\cite{Wa1}): $\Hilb^n(\cit^2//\Gamma) \to \Sym^n(\cit^2/\Gamma),$
where  $\Gamma \subset SL(2)$ is a finite sub-group and  $ \cit^2//\Gamma \to \cit^2/\Gamma $ is the minimal resolution. 
The simplest case is $n=1$.  It can be  shown (\cite{Slo}) that a transverse slice of the sub-regular 
nilpotent orbit in the nilpotent cone 
has $ADE$ singularities, then Grothendieck's simultaneous resolution provides symplectic deformations
of the minimal resolution (see also \cite{GK} section 3).

The next simple case is $n=2$ and $\Gamma = \pm 1$, i.e. the resolution
$\pi: \Hilb^2(T^*\pit^1) \to \Sym^2(\cit^2/{\pm 1})$. Our aim of this note is to construct a symplectic deformation of 
the resolution $\pi$. The key idea is to recover $\pi$ as a slice of some Springer resolution. More precisely,  
let us consider the following two nilpotent orbits in $\mathfrak{sp}_6$:
$$
 \0_{[2,2,2]}: =  \Sp_6 \cdot  \left( \begin{array}{cccccc} 0&  0 & 0 & 1& 0& 0 \\  0&  0 & 0 & 0& 1& 0 \\ 
 0&  0 & 0 & 0& 0& 1 \\ 0 &  0 & 0 & 0& 0 & 0 \\ 0 &  0 & 0 & 0& 0 & 0\\ 0 &  0 & 0 & 0& 0 & 0
\end{array} \right), \quad
\0_{[4,2]} : =  \Sp_6 \cdot \left( \begin{array}{cccccc} 0 &  1 & 1 & 1& 0 & 0 \\0&  0 & 0 & 0& 1& 0 \\ 
 0&  0 & 0 & 0& 0& 1  \\ 0 &  0 & 0 & 0& 0 & 0 \\ 0 &  0 & 0 & -1& 0& 0 \\ 0 &  0 & 0 & -1& 0& 0
 \end{array} \right),
$$
then their closures in $\mathfrak{sp}_6$ are given by:
$$  \overline{\0}_{[2,2,2]}:= \{A \in \mathfrak{sp}_6 \mid A^2 = 0\}, \qquad
\overline{\0}_{[4, 2]} : = \{A \in \mathfrak{sp}_6 \mid A^4 = 0\}. $$

We will prove that the wreath product $\Sym^2(\cit^2/{\pm 1})$ is in fact isomorphic to the intersection of 
a transverse slice of the nilpotent orbit $\0_{[2,2,2]}$
 with the nilpotent orbit closure $\overline{\0}_{[4, 2]} $ 
in $\mathfrak{sp}_6$.
The singular variety $\overline{\0}_{[4, 2]}$ admits exactly two symplectic resolutions. By restricting them to
the transverse slice, we recover exactly the two symplectic resolutions of $\Sym^2(\cit^2/{\pm 1})$. Using
 deformations of Springer resolutions (e.g. \cite{Fu}), we construct a symplectic deformation of $\pi$. 

It is somewhat surprising that we can recover the wreath product $\Sym^2(\cit^2/{\pm 1})$
from nilpotent orbits, although the interplay between nilpotent orbits and Hilbert schemes has
 been noticed  in 
\cite{Man}, where a transverse slice to the nilpotent orbit $\0_{[2m-n,n]} (n \leq m)$
in the nilpotent cone of $\mathfrak{sl}_{2m}$ 
is recovered as an open subset (whose complement is of codimension 1 when $n \geq 2$) of 
the Hilbert scheme $Hilb^n(A_{2m})$, for some
singular surface $A_{2m}$. Here $\0_{[2m-n,n]}$ consists of nilpotent matrices $A \in \mathfrak{sl}_{2m}$
whose Jordan form has only two blocks,  with sizes $2m-n$ and $n$ respectively.
It would be very interesting to  recover other wreath products as a transverse slice
to nilpotent orbits, which would in turn reveal more the mysterious relationships between Hilbert-Chow resolutions
and Springer resolutions, although the two objects are studied in usual separately.

\vspace{0.3 cm}
   
{\bf 1. The transverse slice} \vspace{0.3 cm}

Let $\g$ be a simple Lie algebra and $G$ its adjoint group. 
For any nilpotent element $x \in \g$, by the theorem of Jacobson-Morozov, 
there exists an $\mathfrak{sl}_2$-triplet
$(x, y,h)$. Then $S = x + \g^y$ is a transverse slice to the nilpotent orbit 
$G \cdot x$ in $\g$, and the morphism $G \times S \to \g$ is
smooth (\cite{Slo}, Section 7.4). Here $\g^y := \{ z \in \g \mid [z,y] = 0\}$.

From now on, let $\g = \mathfrak{sp}_6$, and consider the following  
 $\mathfrak{sl}_2$-triplet associated to  the nilpotent orbit $\0_{[2,2,2]}$:
\begin{equation}\label{sl2}
 x_0 = \left( \begin{array}{cc} 0 & I \\ 0 & 0 \end{array} \right), \quad
y_0=\left( \begin{array}{cc} 0 & 0 \\ I & 0 \end{array} \right), \quad
h_0 = \left( \begin{array}{cc} I & 0 \\ 0 & -I \end{array} \right),
\end{equation}
where $I$ is the $3 \times 3$ identity matrix. Note that $\0_{[2,2,2]} = \Sp_6 \cdot x_0.  $

The transverse slice to the orbit $\0_{[2,2,2]}$ is given by
\[
S = x_0 + \g^{y_0} = 
\{ \left( \begin{array}{cc} Z_1 & I \\ Z_2 & Z_1  \end{array}  \right) \mid Z_1 + Z_1^T = 0, Z_2 = Z_2^T \}
\subset \mathfrak{sp}_6.
\]
We choose the following parameters for $Z_1$ and $Z_2$:
\[ 
Z_1 = \left( \begin{array}{ccc} 0 & a_3/2 & -a_2/2 \\ -a_3/2 & 0 & a_1/2 \\
a_2/2 & -a_1/2 & 0 \end{array} \right), \qquad
Z_2 = \left( \begin{array}{ccc} x_1 & y_1 & y_2 \\ y_1 & x_2 & y_3 \\ y_2 & y_3 & x_3
\end{array} \right).
 \]

Note that $\0_{[2,2,2]} \subset \overline{\0}_{[4,2]}$, and the codimension is
is 4. Let $T$ be the scheme intersection $S \cap \overline{\0}_{[4,2]}.$
The variety $ \overline{\0}_{[4,2]}$ is normal and the morphism 
$G \times T \to \overline{\0}_{[4,2]}$ is smooth. It follows that $T$ is normal.
As easily seen, a matrix $A \in S$ is in $T$ if and only if 
$\rk(A) \leq 4$ and $\tr(A) = \tr(A^2) = \tr(A^3) = \tr(A^4) = 0. $

Notice that $\rk(A) \leq 4$ is equivalent to $\rk(Z_2 - Z_1^2) \leq 1.$
The matrix $Z_2 - Z_1^2$ is symmetric, so this is equivalent to the existence
of $u =(u_1, u_2, u_3) \in \cit^3$ such that $Z_2 - Z_1^2 = u^T u,$
from which we can substitute the variables $x_i, y_j$ by $u_k$.
Remark that $u$ and $ -u$ give the same $Z_2$, so we should quotient
 by the following action of $\zit_2$: $u \mapsto -u$.

Now a direct calculus shows that $\tr(A) = \tr(A^3) = 0$.
That $\tr(A^2) = 0$ is equivalent to $ \sum_{i=1}^3 u_i^2 = \sum_{i=1}^3 a_i^2$,
and
$ \tr(A^4) =  2 \tr(Z_1^4) + 2 \tr(Z_2^2) + 12 \tr(Z_1^2 Z_2) = 0$
is equivalent to $\sum_{i=1}^3 a_i u_i = 0$, 
which gives: 
\[
T =  \{(a_1, a_2, a_3, u_1, u_2, u_3) \mid  \sum_i u_i^2 = \sum_i a_i^2, \sum_i a_i u_i = 0 \} /\zit_2,
\] 
where the action of $\zit_2$ is given by $$(a_1, a_2, a_3, u_1, u_2, u_3) \mapsto (a_1, a_2, a_3, -u_1, -u_2, -u_3).$$

Consider the following two nilpotent orbits in $\mathfrak{sp}_6$:
$$ \0_{[3,3]} = \{ A \in \mathfrak{sp}_6 \mid A^3 = 0, \rk(A) = 4\},$$  
$$ \0_{[4,1,1]} = \{ A \in \mathfrak{sp}_6 \mid A^4 = 0, \rk(A) = 3, A^2 \neq 0\}. $$

Then $\0_{[2,2,2]} = \overline{\0}_{[3,3]} \cap \overline{\0}_{[4,1,1]}$ and
 $\overline{\0}_{[3,3]} \subset \overline{\0}_{[4,2]}  \supset \overline{\0}_{[4,1,1]}$. The relationship of inclusions
can be resumed in the following diagram:
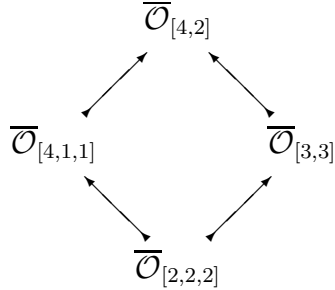
\begin{diagram}
&         &    \overline{\0}_{[4,2]}  &   & \\
& \ruInto &                           &  \luInto   & \\
\overline{\0}_{[4,1,1]}  &   & &  &  \overline{\0}_{[3,3]} \\
& \luInto &                           &  \ruInto   & \\
&         &    \overline{\0}_{[2,2,2]}  &   &
\end{diagram}

 The intersection of $T$ with the two orbit closures $\overline{\0}_{[3,3]}, \overline{\0}_{[4,1,1]}  $
is exactly the singular locus of $T$, which is defined by the following
$$
T \cap \overline{\0}_{[3,3]} = \{ \sum_i a_i^2 =0,  u^T u = -4 Z_1^2   \}, 
$$
$$
T \cap \overline{\0}_{[4,1,1]} = \{ \sum_i a_i^2 = 0,  u_1=u_2=u_3=0\}. 
$$
Both are isomorphic to the surface $\cit^2/{\pm 1}$ with an isolated $A_1$-singularity. The intersection
$T \cap \0_{[2,2,2]} = x_0$ is just a point.
\vspace{0.3 cm}

{\bf 2. The wreath product} \vspace{0.3 cm}

Now we consider the simplest wreath product $\Sym^2 (\cit^2/{\pm 1} ) = \cit^4 / H $,
where $H$ is the subgroup of $\Sp(4)$  generated by the following elements:
$$
\sigma(x_1, x_2, y_1, y_2) = (y_1, y_2, x_1, x_2),  \  \tau(x_1, x_2, y_1, y_2) = (-x_1, -x_2, y_1, y_2).
$$

To write down equations for this affine normal variety, we put
$$
a_1 = - i (x_1^2+y_1^2+x_2^2+y_2^2)/2,  a_2 = (x_1^2+y_1^2-x_2^2-y_2^2)/2, a_3 = x_1x_2 + y_1y_2,
$$
$$
u_1 = x_1y_1 + x_2y_2, u_2 = i(x_1y_1 - x_2y_2), u_3 = i(x_1y_2 + x_2y_1).
$$

The functions $a_1, a_2, a_3$ are $H$-invariant, and the action of $H$ restricts
to a $\zit_2$-action  on $(u_1, u_2, u_3)$ given by $(u_1, u_2, u_3) \mapsto (-u_1, -u_2, -u_3)$.
Now it is straight-ward to check that 
$$ \cit^4 / H \simeq \{(a_1, a_2, a_3, u_1, u_2, u_3) \mid  \sum a_i^2 = \sum u_i^2,\  \sum_i a_i u_i =0 \}/\zit_2,  $$
where the $\zit_2$-action is given by $(a_1, a_2, a_3, u_1, u_2, u_3) \mapsto (a_1, a_2, a_3, -u_1, -u_2, -u_3).  $
 This gives the following proposition.
\begin{Prop}
The transverse slice $T$ is isomorphic to the wreath product $X:= \Sym^2 (\cit^2/{\pm 1} ). $
\end{Prop}

The singular locus of the  wreath product $X$  has two components:
one is the diagonal $\Delta$ and the other will be denoted by
$\Xi$. One sees that the isomorphism between $T$ and $X$ sends $T \cap \overline{\0}_{[3,3]}$ to $\Delta$ and 
$T \cap \overline{\0}_{[4,1,1]}$
to $\Xi$. 

A symplectic resolution of $X$ is given by the composition: 
$$\pi: \Hilb^2(T^* \pit^1)  \to \Sym^2(T^* \pit^1) \to \Sym^2(\cit^2/{\pm 1}) = X.$$
The central fiber of $\pi$ contains a $\pit^2$, so we can blow up this $\pit^2$ and then 
blow down along the other direction, i.e. we can perform a Mukai flop, which gives another
symplectic resolution $\pi^+: \Hilb^2(T^* \pit^1) \to X$.  
One sees that $\pi^{-1}(0) \subset \overline{\pi^{-1}(\Delta - \{0\})}$,
but  $\pi^{-1}(0)$ is not contained in $ \overline{\pi^{-1}(\Xi- \{0\})}$, so $\Delta$ and $\Xi$ are not symmetric
with respect to $\pi$.
For $\pi^+$, it changes the role of 
$\Delta$ and $\Xi$. It is known that any projective symplectic resolutions of $X$ is isomorphic to $\pi$ or $\pi^+$
 (for details, see (\cite{FN}, \cite{Fuj})).
 \vspace{0.3 cm}

{\bf 3. Springer resolutions} \vspace{0.3 cm}

The nilpotent orbit closure $\overline{\0}_{[4,2]}$ admits  exactly two symplectic resolutions, given by
Springer maps: \begin{equation}\label{Springer}
 T^*(G/P_1) \xrightarrow{\phi_1} \overline{\0}_{[4,2]} \xleftarrow{\phi_2} T^*(G/P_2),
\end{equation}
where $P_1$ (resp. $P_2$) is the standard parabolic sub-group of $G$ with flag type 
[1,2,2,1] (resp. [2,1,1,2]). The matrix forms of $P_1$ and $P_2$ are as follows:
$$
P_1 = \{ \left( \begin{array}{cccccc} * &  * & * & *& *& * \\ 0 &  * & * & *& *& * \\ 
0 &  * & * & *& *& * \\ 0 &  0 & 0 & *& 0 & 0 \\ 0 &  0 & 0 & *& *& * \\ 0 &  0 & 0 & *& *& *
\end{array} \right) \}, \quad
P_2= \{\left( \begin{array}{cccccc} * &  * & * & *& *& * \\ * &  * & * & *& *& * \\ 
0 &  0 & * & *& *& * \\ 0 &  0 & 0 & *& * & * \\ 0 &  0 & 0 & *& *& * \\ 0 &  0 & 0 & 0& 0& *
 \end{array} \right) \}.
$$

The restrictions of $\phi_1, \phi_2$ to the pre-image of the transverse slice $T$ give two
 projective  symplectic resolutions of $T$.
$$ \Hilb^2(T^*\pit^1) \simeq Z_1 \xrightarrow{\pi_1} T \xleftarrow{\pi_2} Z_2 \simeq \Hilb^2(T^* \pit^1) .$$
\begin{Prop}
The two symplectic resolutions $\pi_1, \pi_2$  are related by a Mukai flop, in particular, they are not isomorphic.
Furthermore $\pi_1 = \pi$ and $\pi_2 = \pi^+$.
\end{Prop}
\begin{proof}
We will calculate the central fiber over the point $x_0 =  T \cap \0_{[2,2,2]}$ (c.f. (\ref{sl2}))
 under the maps $\pi_1$ and $\pi_2$.
Let $\{e_i, 1 \leq i \leq 6  \}$ be the natural basis of $\cit^6$ and the symplectic form
is $\omega = \sum_{i=1}^3 e_i^* \wedge e_{i+3}^*.$
Notice that $Im(x_0) = Ker(x_0) = \cit \langle e_1, e_2, e_3 \rangle   =: K$ is Lagrangian.

It is easy to see that $$ \pi_1^{-1}(x_0) = \{ \text{flags} \ (F_1 \subset F_2)  \mid  x_0 F_2 \subset F_1 \subset
K, F_2 = F_2^\perp, \dim F_1 = 1 \}.$$
 Since $x_0 F_2 \subset F_1$ is of dimension 1, one has two possibilities:

(i). $ dim (K \cap F_2) =  2$, then $F_1 = x_0 F_2 \subset x_0 F_1^\perp$. Suppose that $F_1$ is generated by
$\sum_{i=1}^3 a_i e_i$, then $x_0 F_1^\perp = \{ \sum_{i=1}^3 b_i e_i| \sum_i a_i b_i = 0 \}.$ The condition 
$F_1  \subset x_0 F_1^\perp$ is equivalent to $\sum_{i=1}^3 a_i^2 = 0$, which is a $\pit^1$ inside $\pit(K)$.
The condition for $F_2$ is just $x_0 F_1^\perp \subset F_2  \subset x_0^{-1} F_1$ which is a $\pit^1$. 
So finally this component is a $\pit^1$-bundle over $\pit^1$.

(ii). $F_2 = K$, then this is isomorphic to $\pit(K)$. The two components intersect at a curve $C_1 \simeq \pit^1$
 inside $\pit(K)$.

 The fiber of $x_0$  under $\pi_2$
consists of flags $(F_1 \subset F_2)$ such that $x_0 F_2 \subset F_1 \subset
K,  F_2 = F_2^\perp$ and $\dim F_1 = 2$. Since $F_1 \subset K \cap F_2$, so $\dim
(K \cap F_2) \geq 2$. There are two cases:

(i).  $dim (K \cap F_2) =  2$, then $F_1 = K
\cap F_2$ and $x_0 F_2 \subset F_1.$ This gives that $x_0 F_2 = x_0 F_1^\perp \subset F_1$.
 Suppose $F_1$ is generated by $\sum_{i=1}^3 a_i e_i, \sum_{i=1}^3 b_i e_i$.
Then we have $x_0 F_1^\perp = \{\sum_i c_i e_i | \sum_i a_i c_i = \sum_i b_i c_i = 0 \}.$
The condition $x_0 F_1^\perp \subset F_1$ is equivalent to the existence of $(y, y') \neq (0,0)$
such that $y (\sum_i a_i^2) + y' (\sum_i a_i b_i) = 0  $ and $y (\sum_i a_i b_i) + y' (\sum_i b_i^2)=0$.
So the condition for $F_1$ is $(\sum_i a_i^2) (\sum_i b_i^2) = (\sum_i a_i b_i)^2$.
Under the Pl\"ucker embedding $\pit(\wedge^2 F_1) \to \pit(\wedge^2 K) \simeq \pit(K^*)$, one sees that this is
a conic in $\pit^2$.   The condition for $F_2$ turns to be $F_1 \subset F_2 \subset F_1^\perp$. 
So this component is a $\pit^1$-bundle over $\pit^1$.

(ii). $K= F_2$, then $F_1 \subset K$, this component is just $\pit(K^*)$. The two components intersects
at a $C_2 \simeq \pit^1$ inside $\pit(K^*)$.

Now it is clear that the two resolutions are different and are  related by the Mukai flop along the component $\pit(K^*)$, and
$C_1, C_2$ are dual conics.

Now we will identify $\pi_1, \pi_2$ with $\pi, \pi^+$. By definition, we have
$$\pi_1^{-1}(T \cap \overline{\0}_{[3,3]}) = \{(F_1 \subset F_2, z) \mid z F_2 \subset F_1 \subset Ker(z), F_2 = F_2^\perp, \dim F_1 = 1\},$$
where $z$ is in $T \cap \overline{\0}_{[3,3]}$.
Consider the elements  $z_t \in T \cap \overline{\0}_{[3,3]}$ $(t \in \cit)$ given by 
\[
z_t= \left( \begin{array}{cc} t B & I \\ -3t^2 B^2 & t B  \end{array}  \right), \text{with} \quad
B = \left( \begin{array}{ccc} 0 & \sqrt{-1} & 1 \\ - \sqrt{-1} & 0 & 0 \\ -1 & 0 & 0 \end{array}  \right).
\]
One has $Ker(z_t) = \cit \langle  e_1 + \sqrt{-1}t e_5 + t e_6, e_2 - \sqrt{-1} e_3 \rangle$.
 When
$t$ goes to $0$, $Ker(z_t)$ goes to $\cit \langle e_1, e_2 - \sqrt{-1} e_3\rangle,$  
thus the limit of $\pi_1^{-1}(z_t)$
will be $\pit(\cit \langle e_1, e_2- \sqrt{-1} e_3 \rangle ) \subset \pit(K)$, which is not a point. 
This shows that $\pi_1= \pi$ by the description of $\pi$ and $\pi^+$ in section 2.
\end{proof}
\vspace{0.3 cm}

{\bf 4. Symplectic deformations} \vspace{0.3 cm}

A deformation of the symplectic resolutions $\phi_i, i=1, 2$ (cf. \eqref{Springer}) can be constructed as follows(\cite{Fu}). 
Let $\mathfrak{c}_i$ be the center of the Levi sub-algebra of $\p_i: = Lie(P_i)$ and $\u_i$
the nil-radical of $\p_i$. The vector space $V_i: =  \mathfrak{c}_i + \u_i$ is a flat family over $\mathfrak{c}_i$.
Let $\Y_i$ be the closed sub-variety 
$$\Y_i: = \{ (z, v) \in \mathfrak{c}_i \times G\cdot V_i \mid v \in G \cdot (z + \u_i) \}. $$  
Now consider the morphism $\Phi_i: G \times^{P_i} V_i \to \Y_i$ given by $$g * (z + u) \mapsto (z, g \cdot (z+u)),$$
where $g\in G, z\in \c_i$ and $u \in \u_i$.
Notice that if $z \neq 0$, then  $z + \u_i = P_i \cdot z$, so this morphism is well-defined. 
One can show that $\Phi_i$ is birational and
it gives a family of morphisms over $\mathfrak{c}_i$.
When $z \in \mathfrak{c}_i$ is generic, in the sense that the stabilizer $G^z$ of $z$ is exactly
 the Levi sub-group $L_i$  of $P_i$, then
$$\Phi_i^z: G \times^{P_i} (z + \u_i) \simeq G \times^{P_i} (P_i \cdot z) \to \Y_i^z = G\cdot z \simeq G/L_i$$
is an isomorphism. Notice that $\Y_i^z$ is a semi-simple orbit, thus it is symplectic.
 When $z = 0$, the map $\Phi_i^0$ is just
the Springer resolution $ \phi_i. $ In other words, $\Phi_i$ is a symplectic deformation of $\phi_i$
with base $\c_i$.

Let $\T_i$ be the intersection $ (\mathfrak{c}_i  \times S) \cap \Y_i$ and $\Z_i$ its pre-image under the morphism
$\Phi_i$, which gives a map $\Z_i \xrightarrow{\Pi_i} \T_i$ over $\mathfrak{c}_i$. 
Now we will show that the family $\psi_i: \Z_i \to \mathfrak{c}_i$ is smooth. Recall that the map $G \times S \to \g$
is smooth, so is $G \times (S \cap (G \cdot V_i)) \to G \cdot V_i$ (\cite{Slo} Section 5.1). The morphism 
$\Phi_i$ is $G$-equivariant, so $G \times \Z_i \to G \times^{P_i} V_i$ is smooth. Notice that the map
$ G \times^{P_i} V_i \to \mathfrak{c}_i $ is smooth, so is the composition $G \times \Z_i \to \mathfrak{c}_i$.
The projection $G \times \Z_i \to \Z_i$ is smooth, which implies the smoothness of $\Z_i \to \mathfrak{c}_i$.

An immediately corollary is that $\Z_i$ is smooth and 
 for any $0 \neq z \in \mathfrak{c}_i$ generic, the intersection $S \cap G\cdot z$ 
is smooth and symplectic, which deforms $\Sym^2(\cit^2/{\pm 1})$.
So $\Z_i \xrightarrow{\Pi_i} \T_i$ gives a symplectic deformation
of the symplectic resolution $\pi_i: \Hilb^2(T^* \pit^1) \to \Sym^2(\cit^2/{\pm 1}).$ 

\vspace{0.3 cm}

{\bf 5. Universal Poisson deformations} \vspace{0.3 cm}

Now we will show that our picture is similar to that of Brieskorn (see section 3 \cite{GK}).
We will only consider $\phi_1$. 
 To simplify the notations, we will write $P$ (resp. $\phi$, $L$ etc.)
instead of $P_1$ (resp. $\phi_1$, $L_1$ etc.).

 Fix a maximal torus $U$ in $G$ and a Cartan sub-algebra $\h$ in $\g=\mathfrak{sp}_6$. 
Coordinates in $\h$ are denoted by $(h_1, h_2, h_3)$.
We define the Weyl group of $L$ to be  $W(L): = N_L(U)/U$, where $N_L(U)$ is the normalizer of 
$U$ in $L$. The partial Weyl group of $P$ is $W^P: = N_G(L)/L$. Then $W^P$ is naturally isomorphic
to the quotient $ N_{W(G)}(W(L))/W(L),$ where $W(G)$ is the Weyl group of $G$.  

It is easy to see that $W(L)$ is isomorphic to $\zit_2$, acting on $\h$ by
$(h_1, h_2, h_3) \mapsto (h_1, h_3, h_2)$.
The center $\c$ of $Lie(L)$ is naturally identified with the fixed point set $\h^{W(L)}.$
The group $W^P$ is isomorphic to $\zit_2 \times \zit_2$, which acts on $ \h^{W(L)} $
by $(h_1, h_2, h_2) \mapsto (-h_1, h_2, h_2)$ and $(h_1, h_2, h_2) \mapsto (h_1, -h_2, -h_2),$
i.e. it is the sum of two copies of  the sign representation of $\zit_2$.

Let $\T'$ be the intersection $S \cap G \cdot (\c + \u)$, then we have a natural projection
$p: \T \to \T'  $ and the following diagram is commutative:
\begin{equation}\label{eq2}
\begin{CD}
\Z @>\Pi>> \T @>p>> \T' \\
@V\psi VV    @V\psi' VV @V\beta VV \\
\h^{W(L)} @>id>> \h^{W(L)}  @>\eta >> \h^{W(L)}/W^P,
\end{CD}
\end{equation}
where $\eta$ is the natural quotient map and $\beta$ is the restriction to $\T'$ of the adjoint 
quotient map $ \g \to \h/W(G)$. \vspace{0.2 cm}

{\em Claim:} The second Poisson cohomoloy $HP^2(T)$ can be naturally identified with $\h^{W(L)}/W^P$. 
\vspace{0.2 cm} 

Let $H'$ be the semi-direct product of $\zit_2$ with $W^P$. Let $\zit_2$ acts on $\h^{W(L)}$ by 
$(h_1, h_2, h_2) \mapsto (h_2, h_1, h_1)$, then it is easy to see that $(\h^{W(L)} \oplus (\h^{W(L)})^*)/H'$
is isomorphic to $T \simeq \Sym^2(\cit^2/{\pm 1}).$ By \cite{GK} (section 4), we have
$HP^2(T)$ is naturally isomorphic to 
$HP^2(T \cap \overline{\0}_{[3,3]}) \oplus HP^2(T \cap \overline{\0}_{[4,1,1]})$.
By Lemma 3.1 \cite{GK}, $HP^2(T \cap \overline{\0}_{[3,3]})$ is naturally identified with 
$ \cit v_1 / \zit_2$ and $HP^2(T \cap \overline{\0}_{[4,1,1]})$ is identified with 
$\cit v_2 / \zit_2$, where $v_1 = (0, 1, 1), v_2 = (1, 0, 0)$ are two points in $\h^{W(L)}$,
and the group $\zit_2$ acts by sign representation.

Note that under this identification, $\h^{W(L)}$ is  identified with 
$H^2(\Hilb^2(T^* \pit^1))$.
The first square in \eqref{eq2} is the symplectic deformation
of $\pi$, the second square is Cartesian. 
For the vertical morphisms, $\psi$ is a universal Poisson deformation of $\Hilb^2(T^*\pit^1)$, 
$\psi'$ is similar to the Calogero-Moser  deformation
of $T \simeq (\h^{W(L)} \oplus (\h^{W(L)})^*)/H'$, 
 and $\beta$ is the universal Poisson deformation of $T$. 
\begin{Rque}
A diagram analogous to \eqref{eq2} can be
constructed using the same method for any Springer resolution.
\end{Rque}
\vspace{0.3 cm}

{\em Acknowledgments:}  This note comes out from very helpful discussions with S.-S. Roan. I want to 
thank him and W. Wang for their interest in my work, and M.S.R.I.(Berkeley) and Max-Planck-Institut f\"ur Mathematik (Bonn)
 for the hospitality.

\quad \\
C.N.R.S.,
Laboratoire J. Leray (Math\'ematiques)\\
 Facult\'e  des sciences, Univ. de Nantes \\
2, Rue de la Houssini\`ere,  BP 92208 \\
F-44322 Nantes Cedex 03 - France\\
\quad \\
fu@math.univ-nantes.fr
\end{document}